\documentclass[a4paper,12pt,leqno]{article}
\usepackage[utf8]{inputenc}
\usepackage[T1]{fontenc}
\usepackage[french, english]{babel}
\usepackage{typearea}
\usepackage{color}
\usepackage{enumerate, enumitem, moreenum}

\usepackage{amstext,amsthm,a4,amscd,amssymb}
\usepackage{amsmath,amsfonts,}

\usepackage{textcomp}
\usepackage{tikz}
\usepackage{fancyhdr}
\usepackage{url}
\usepackage{mathrsfs}
\usepackage{dsfont}
\usepackage{cancel}
\usepackage{ulem}
\usepackage[justification=centering]{caption}

\usepackage
{hyperref}
\hypersetup{
    colorlinks = 
    true,
    linkcolor={black},
    linkbordercolor = {white},
    citecolor={black},
     }

\usepackage[all]{xy}

\usepackage[a4paper,top=2cm,bottom=2.5cm,left=2cm,
right=2cm,footskip=1cm]
{geometry}

\usepackage{tikz}
 \usetikzlibrary{patterns}

\numberwithin{equation}{section}

\newtheorem{df}{Definition}[section]

\newtheorem{thm}[df]{Theorem}

\newtheorem{lem}[df]{Lemma}
\newtheorem{crl}[df]{Corollary}
\newtheorem{rmk}[df]{Remark}

\newcommand{\C}{\mathbb{C}}

\newcommand{\N}{\mathbb{N}}
\newcommand{\D}{\mathbb{D}}

\newcommand{\bL}{\boldsymbol{L}}

\newcommand*{\LargerCdot}{\raisebox{-0.25ex}
{\scalebox{1.5}{$\cdot$}}}

\newcommand{\comment}[1]{}

\title{Fubini\textendash Study forms on punctured Riemann surfaces \bigbreak
\Large \textit{Formes de Fubini\textendash Study sur des surfaces de Riemann épointées}}
\author{\normalsize\textsc{
Razvan Apredoaei
\footnote{ Universit\'e Paris Cit\'e, CNRS, IMJ-PRG,
B\^atiment Sophie Germain, UFR de Math\'ematiques,
Case 7012, 75205 Paris Cedex 13, France.  Email: apredoaei@imj-prg.fr }
,
Xiaonan Ma\footnote{ Chern Institute of Mathematics and LPMC,
Nankai University, Tianjin 300071, People's Republic of China.
Email: xiaonan.ma@nankai.edu.cn }
,
Lei Wang\footnote{ School of Mathematics and Statistics,
Huazhong University of Science and Technology,
Wuhan, Hubei 430074, P.R.China.  Email: wanglei2017@hust.edu.cn}}}
\date{}

\comment{
\author{Razvan Apredoaei}
\address{Universit\'e Paris Cit\'e, CNRS,
IMJ-PRG, B\^atiment Sophie Germain, UFR de Math\'ematiques,
Case 7012, 75205 Paris Cedex 13, France}
\email{apredoaei@imj-prg.fr}

\author{Xiaonan Ma}
\address{Chern Institute of Mathematics and LPMC,
Nankai University, Tianjin 300071, People's Republic of China}
\email{xiaonan.ma@nankai.edu.cn}

\thanks{Ma is partially supported by NSFC  No.11829102,
ANR-21-CE40-0016 and funded through the Institutional Strategy of
the University of Cologne within the German Excellence Initiative.}

\author{Lei Wang}
\address{School of Mathematics and Statistics,
Huazhong University of Science and Technology, Wuhan,
Hubei 430074, P.R.China}
\email{wanglei2017@hust.edu.cn}
}

\makeatletter
\def\blfootnote{\xdef\@thefnmark{}\@footnotetext}
\makeatother

\begin{document}
\blfootnote{
X. M. is partially supported by Nankai Zhide Foundation, 
ANR-21-CE40-0016 and funded through the Institutional Strategy of
the University of Cologne
within the German Excellence Initiative.
L. W. is partially
supported by NSFC No. 11801187.
           \vspace{4pt}
         \footnoterule
           \vspace{1pt}}

\makeatletter
\renewcommand\section{\@startsection {section}{1}{\z@}%
                           {-3.5ex \@plus -1ex \@minus -.2ex}%
                          {2.3ex \@plus.2ex}%
                      {\centering\sc\normalsize}}
\renewcommand\subsection{\@startsection{subsection}{2}{\z@}%
                           {-3.25ex\@plus -1ex \@minus -.2ex}%
                           {1.5ex \@plus .2ex}%
                          {\normalsize\sf}}
\renewcommand\subsubsection{\@startsection{subsubsection}{3}
{\z@}%
                     {-3.25ex\@plus -1ex \@minus -.2ex}%
                    {1.5ex \@plus .2ex}%
                   {\normalsize\it}}

\makeatother

\maketitle

\selectlanguage{english}
\begin{abstract}
In this paper we consider a punctured Riemann surface
endowed with a Hermitian metric that equals the Poincar\'e
metric near the punctures, and a holomorphic line bundle
that polarizes the metric. We show that the quotient of the induced Fubini\textendash Study forms by Kodaira maps
of high tensor powers of the line bundle and the Poincar\'e
form near the singularity
 grows polynomially uniformly on a neighborhood of the
 singularity as the tensor power tends to infinity,
 as an application of the method in \cite{AMM22}.
\end{abstract}

\selectlanguage{french}
\begin{abstract}
Dans cet article, nous considérons une surface de Riemann 
épointée munie d'une métrique hermitienne qui coïncide 
avec la métrique de Poincaré près des points de ponction, 
ainsi qu'un fibré en droites holomorphe qui polarise la métrique. Nous montrons que le quotient des formes induites de Fubini\textendash Study 
par les applications de Kodaira des puissances tensorielles élevées 
du fibré en droites et de la forme de Poincaré près de la singularité
croît de manière polynomiale et uniforme dans un voisinage 
de la singularité lorsque la puissance tensorielle tend vers l'infini, 
en application de la méthode décrite dans \cite{AMM22}.
\end{abstract}


\selectlanguage{english}


\section{Introduction}
In this paper we study the asymptotics of the induced
Fubini\textendash Study metrics by Kodaira maps of high tensor powers
of a singular Hermitian line bundle over a Riemann surface
under the assumption that the curvature has singularities
of Poincar\'e type at a finite set.
We show namely that the {\it quotient} of the induced
Fubini\textendash Study metrics by Kodaira maps of high tensor powers
of the line bundle and the Poincar\'e model near
the singularity grows polynomially uniformly on
a neighborhood of the singularity as the tensor power
tends to infinity.
In \cite{AMM16,AMM21}, Auvray, Ma and Marinescu obtained
a weighted estimate in the $C^m$-norm near the punctures
for the difference of the global Bergman kernel and of the
 Bergman kernel of the Poincar\'e model near the singularity,
 uniformly in the tensor power $p$ of the given line bundle;
and an application in arithmetic geometry 
is the sharp uniform asymptotics 
of the sup-norm of automorphic (cusp) forms 
associated with a Fuchsian group of the first
kind which was studied by Ullmo, Kramer $\dots$
\cite{MiUll98, JK04, fjk}.
Moreover, in \cite{AMM22} they also show that the quotient
of the global Bergman kernel and the Bergman kernel of
the Poincar\'e model near the singularity tends to one
up to $\mathcal{O}(p^{-\infty})$.
This paper is an application of the results and the method
of \cite{AMM21, AMM22}.

Following \cite{Tian}, an important application of the
expansion of the Bergman kernel is the convergence of
the induced Fubini\textendash Study forms by Kodaira maps.
For more works on Bergman kernels,
cf. \cite{Catlin, DLM06, MM08, Zeld, Finski24}, and  we refer the readers
  to the book \cite{MM07} for a comprehensive study of
  Bergman kernels.
Donaldson \cite{Donaldson} uses the expansion of
the Bergman kernel to study the relation between constant
scalar curvature K\"ahler metrics and Chow stability.
Coming to our context, an interesting problem is the
relation between the existence of special complete/singular
metrics and the stability of the pair $(X,D)$ where
$D$ is a smooth divisor of a compact K\"ahler manifold $X$;
 see e.g. the suggestions of \cite[Section 3.1.2]{Szeke}
 and \cite{LW},  for the case of ``asymptotically hyperbolic
  K\"ahler metrics'', which naturally generalize to higher
  dimensions the complete metrics $\omega$ studied here.

\vskip 10pt
We place ourselves in the setting of \cite{AMM21}
which we describe now.
Let $\overline\Sigma$ be a compact Riemann surface and let
$D=\{a_1,\ldots,a_N\}\subset\overline\Sigma$ be a finite set.
We consider the punctured Riemann surface
$\Sigma = \overline{\Sigma}\smallsetminus D$ and a Hermitian
 form
$\omega_{\Sigma}$ on $\Sigma$.
Let $L$ be a holomorphic line bundle on $\overline{\Sigma}$,
and let $h$ be a singular Hermitian metric on $L$ such that:

   \begin{enumerate}[label = (\greek*)]
    \item\label{alpha} $h$ is smooth over $\Sigma$,
and for all $j=1,\ldots,N$, there is a trivialization of $L$
 in the complex neighborhood $\overline{V_j}$ of $a_j$ in
 $\overline{\Sigma}$, with associated coordinate $z_j$
 such that $|1|_{h}^2(z_{j})= \big|\!\log(|z_j|^2)\big|$;
    \item\label{beta} there exists $\varepsilon>0$ such that the
(smooth) curvature $R^L$ of $h$ satisfies
$iR^L\geq\varepsilon\omega_{\Sigma}$ over $\Sigma$
and moreover, $iR^L=\omega_{\Sigma}$
on $V_j:=\overline{V_j}\smallsetminus\{a_j\}$;
in particular, $\omega_{\Sigma} = \omega_{\D^*}$
in the local coordinate $z_j$ on $V_j$
and $(\Sigma, \omega_{\Sigma})$ is complete.
   \end{enumerate}
Here $\omega_{\D^*}$ denotes the Poincar\'{e} metric on the
punctured unit disc $\D^* = \D \setminus\{0\}$, where $\D$
is the unit disc, normalized as follows:
\begin{equation}   \label{eq1}
\omega_{\D^*} := \frac{idz\wedge d\overline{z}}
{|z|^2\log^2(|z|^2)}\,\cdot
\end{equation}
For $p\geq1$, let $h^p:=h^{\otimes p}$ be the metric induced
 by $h$ on $L^p\vert_{\Sigma}$,
where $L^p:=L^{\otimes p}$. We denote by
$H^0_{(2)}(\Sigma,L^p)$
the space of ${\bL}^2$-holomorphic sections of $L^p$
relative to the metrics $h^p$ and $\omega_\Sigma$
endowed with the obvious inner product.
By \cite[(6.2.17)]{MM07}, the sections from
$H^0_{(2)}(\Sigma, L^p)$ extend to holomorphic sections of
$L^p$ over $\overline\Sigma$ which vanish on $D$.
In particular, the dimension $d_p$ of $H^0_{(2)}(\Sigma, L^p)$
is finite.


We denote by $B_p(\LargerCdot)$ the Bergman kernel function
of the orthogonal projection $B_{p}$ from the space
of $\bL^{2}$-sections of $L^{p}$ over $\Sigma$ onto
$H^0_{(2)}(\Sigma,L^p)$.
If $\{S_\ell^p\}_{\ell=1}^{d_p}$ is an orthonormal
basis of $H^0_{(2)}(\Sigma,L^p)$, then
\begin{equation}\label{eq4}
B_p(x):=\sum_{\ell=1}^{d_p}|S^p_\ell(x)|_{h^p}^2\,.
\end{equation}
Note that these are independent of the choice of basis (see
\cite[(6.1.10)]{MM07}).
Similarly, let $B_p^{\D^*}(x)$ be the Bergman kernel function of
$\big(\D^*, \omega_{\D^*},
\C,\big|\!\log(|z|^2)\big|^p\, h_{0})$ 
with $h_{0}$ the flat Hermitian metric on the trivial line
 bundle $\C$.

We fix a point $\mathbf{a}\in D$ and work in coordinates
centered at $\mathbf{a}$.
Let $\mathfrak{e}_L$ be the holomorphic frame of $L$ near
 $\mathbf{a}$ corresponding to the trivialization in the
 condition \ref{alpha}.
By the assumptions \ref{alpha} and \ref{beta} we have the
following identification of the geometric data in the coordinate
 $z$ on the punctured disc $\mathbb{D}^*_{4r}$ of radius $4r$
 centered at $\mathbf{a}$, via the trivialization
 $\mathfrak{e}_L$ of $L$,
\begin{equation}\label{eq8}
\big(  \Sigma, \omega_{\Sigma}, L,h \big)
 \big|_{\mathbb{D}^*_{4r} }
= \big( \mathbb{D}^*, \omega_{\mathbb{D}^*},
\mathbb{C}, h_{\mathbb{D}^*} = \big| \log(|z|^2) \big|
\cdot h_0 \big)  \big|_{\mathbb{D}^*_{4r} }, \;\;\;
\text{with} \; \; 0< r < (4e)^{-1}.
\end{equation}
In \cite{AMM22},
Auvray, Ma and Marinescu proved the following estimates.

\begin{thm}\cite[Theorems 1.2 and 1.3]{AMM22}
  \label{thm_prevThm}
If $(\Sigma, \omega_\Sigma, L, h)$ fulfill conditions
 \ref{alpha} and \ref{beta},
then for all $k\in \mathbb{N}$, $l\in \mathbb{N}^*$ and
$D_1,\cdots, D_k \in \Big\{ \frac{\partial}{\partial z},
\frac{\partial}{\partial \bar{z}} \Big\}$,
there exists $C>0$ such that for any
$p\in\N^{*}$ we have
\begin{equation}\label{eq1.7}
\sup_{z\in V_1\cup\ldots\cup V_N}
\left|(D_1 \cdots D_k) \bigg( \frac{B_p}{B_p^{\D^*}}(z) - 1 \bigg)
\right| \leq C p^{-\ell}.
\end{equation}
\end{thm}

\vskip 10pt

Let us consider the Kodaira map at level $p\geq 2$ induced
by $H^0_{(2)}(\Sigma, L^p)$, which is a meromorphic map
defined by
\begin{equation}\label{eq5}
J_{p,(2)}: \Sigma \dashrightarrow
\mathbb{P}(H^0_{(2)}(\Sigma, L^p)^*)
 \cong \mathbb{CP}^{d_p-1}, \;\;\; x\mapsto
 \{ \sigma \in H^0_{(2)}(\Sigma, L^p): \sigma(x) = 0\}.
\end{equation}
By \cite[p2364]{AMM22}, $J_{p,(2)}$ is an embedding for $p\gg 1$. 

The ${\bL}^2$-metric on $H^0_{(2)}(\Sigma, L^p)$ induces
 a Fubini\textendash Study form $\omega_{FS, p}$ on the
 projective space
$\mathbb{P}(H^0_{(2)}(\Sigma, L^p)^*) $.
We have by  \cite[Theorem 5.1.3, (5.1.21)]{MM07}
\begin{equation}\label{eq6}
\frac{1}{p} J^*_{p,(2)} \omega_{FS,p} = \frac{i}{2\pi} R^L
+ \frac{i}{2\pi p} \partial \bar{\partial} \log (B_p).
\end{equation}
From (\ref{eq1.7}) and (\ref{eq6})
 (cf. \cite[Theorem 4.1]{AMM22}),  we have, as
 $p\rightarrow +\infty$,
\begin{equation}\label{eq7}
\frac{1}{p} J^*_{p,(2)} \omega_{FS,p}
= \frac{1}{2\pi} \omega_{\Sigma}
 + \frac{i}{2\pi p} \partial \bar{\partial}
 \log (B_p^{\mathbb{D}^* }) + \mathcal{O}(p^{-\infty}),
\end{equation}
uniformly on $V_1 \cup V_2 \cup \cdots \cup V_N$.
Note that, by condition \ref{beta}, we have
\begin{equation}\label{eq001}
  \omega_{\D^*} = \omega_{\Sigma} = i R^L
  = - i\partial\bar{\partial} \log \big| \log(|z|^2) \big| \;\;\;
  \text{on} \;\; V_j.
\end{equation}

\vskip 10pt
The main result of the present paper is the following estimate
of the quotient of the induced Fubini\textendash Study forms
 $\frac{1}{p} J^*_{p,(2)} \omega_{FS,p}$ by Kodaira maps
 of $L^p$ and the Poincar\'e form
 $\omega_{\mathbb{D}^*}(z)$:

\begin{thm} \label{thm_MainThm}
If $(\Sigma, \omega_{\Sigma}, L, h)$ fulfill conditions
\ref{alpha} and \ref{beta}, then
\begin{equation}\label{eq9}
\sup_{z\in V_1 \cup \cdots \cup V_N }
\Bigg | \frac{  J^*_{p,(2)} \omega_{FS,p}(z) }
{p \; \omega_{\mathbb{D}^*} (z)  } \Bigg |
 = \mathcal{O}(p^3) \qquad \text{as }p\to +\infty.
\end{equation}
\end{thm}

Let us mention that the difficulty of the estimate (\ref{eq9})
consists in the fact that $B_p^{\mathbb{D}^*}(\cdot)$ vanishes 
at $0$.
Moreover, by \cite[Theorem 6.1.1]{MM07} and (\ref{eq6}), for any
$K \subset \Sigma$ compact,
\begin{equation}\label{eq0006}
\frac{1}{p} J^*_{p,(2)} \omega_{FS,p} - \frac{i}{2\pi} R^L
= \mathcal{O}(p^{-1}) \;\;\;\;\; \text{on} \;\; K.
\end{equation}

\vskip 10pt
We give an important example where Theorem \ref{thm_MainThm}
applies.
Let $\Gamma \subset PSL(2, \mathbb{R})$ be a geometrically
finite Fuchsian group of the first kind without elliptic
 elements. Then $\Sigma: = \Gamma \backslash \mathbb{H}$
 can be compactified by finitely many points
 $D= \{a_1, \cdots, a_N\}$ into a compact Riemann surface
 $\overline{\Sigma}$ \cite[p956]{AMM21}.
Let $\omega_{\Sigma}$ be the K\"ahler\textendash Einstein metric on
$\Sigma$ induced by the Poincar\'e metric
$\omega_{\mathbb{H}} = \frac{ idz \wedge d\bar{z}}{ 4|{\rm Im} z|^2}$
on the upper-half plane $\mathbb{H}$.

Let $\mathcal{S}^{\Gamma}_{2p}$ be the space of cusp forms
(Spitzenformen) of weight $2p$ of $\Gamma$ endowed with
the Petersson inner product (cf. \cite[p995]{AMM21}).
We can form the Bergman kernel function of
$\mathcal{S}_{2p}^{\Gamma}$ as in (\ref{eq4}),
denoted by $S_{p}^{\Gamma}$\cite[p995]{AMM21},
then
\begin{equation}\label{eq10a}
  S_{p}^{\Gamma}(z) = \sum_j |f_j(z)|^2 (2 \text{Im}z)^{2p}
\end{equation}
with $\{f_j\}$ any orthonormal basis of $\mathcal{S}^{\Gamma}_{2p}$.
For any $p>0$ and $z\in \Sigma$, define
\begin{equation}\label{eq.berg}
  \omega_{\Sigma}^{Ber, p} : =
  \frac{i}{2\pi} \partial \bar{\partial}
  \log (S_{p}^{\Gamma}(z)) .
\end{equation}
As a corollary of Theorem \ref{thm_MainThm}, we get the
following:

\begin{crl}\label{cor}
Let $\Gamma \subset PSL(2, \mathbb{R})$ be a geometrically
finite Fuchsian group of the first kind without elliptic
elements and $\Sigma : =\Gamma \backslash \mathbb{H}$.
Then
\begin{equation}\label{eq10}
\sup_{z\in \Sigma}  \left| \frac{ \omega_{\Sigma}^{Ber, p}(z) }
{ p \; \omega_{\Sigma}(z)} \right|
= \mathcal{O}(p^3) \quad \text{as} \;\; p\rightarrow +\infty.
\end{equation}
\end{crl}

\begin{rmk}
Corollary \ref{cor} still holds if $\Gamma$ has elliptic elements.
In this case, the quotient $\Sigma = \Gamma \backslash \mathbb{H}$
is an orbifold, by the result of Dai\textendash Liu\textendash Ma \cite[(5.25)]{DLM06},
Ma\textendash Marinescu \cite[Section 5.4.3]{MM07} on the Bergman kernel on orbifolds
and the localization of the asymptotics of the Bergman kernel on compact
part of $\Sigma$, from \cite[(1.15), p998]{AMM21}, we know near
the orbifold points, $\mathcal{O}(p^3)$ in (\ref{eq10}) can
be estimated as $\mathcal{O}(1)$.

\end{rmk}

\begin{rmk}
In \cite[Main Theorem 1]{Arya23}, the authors claim a stronger estimate
than (\ref{eq10}) and they prove their
estimate by using a computation on classical modular forms,
 but unfortunately their proof is not complete.
As a comparison, the difficulty in proving (\ref{eq10})
arises from the fact that the Bergman kernel function vanishes on $D$,
but \cite[Proposition 2.3]{Arya23} claims that it has
a uniform strictly positive lower bound,
cf. \cite[Corollary 3.6]{AMM21} and the explanation after
(\ref{eq16}).

In \cite[Main Theorem 2]{Arya24}, the authors claim
a high dimension version of (\ref{eq10}) for
ball quotients, but unfortunately their proof is also not complete.
The similar incomplete argument occurs in the proof of
\cite[Proposition 3.6]{Arya24}. 
\end{rmk}

\section{Proof of Theorem \ref{thm_MainThm} and Corollary \ref{cor}}
In this section,  we establish Theorem \ref{thm_MainThm} and
Corollary \ref{cor}  by using the techniques in \cite{AMM22}.
Our starting point is (\ref{eq7}), and thus we only need to work on
 the punctured unit disc $\mathbb{D}^*$.
The key part in our proof of Theorem \ref{thm_MainThm}
in Section \ref{s2.1}
is Lemma \ref{lem} which is established in Section \ref{s2.2}:
For the four terms in (\ref{eq26}) and (\ref{eq27}),
we get $\mathcal{O}(p^{-\infty})$.
The main estimate $\mathcal{O}(p^3)$
is from (\ref{eq28}) as $B_p^{\D^*}(z)$ oscillates
 too violently near $|z|=e^{-p}$
  (cf.  \cite[Section 3.2]{AMM21}).
In Section \ref{s2.3}, we explain the proof of Corollary \ref{cor}.


\subsection{Proof of Theorem \ref{thm_MainThm}} \label{s2.1}
From \cite[(2.6) and (2.7)]{AMM22}, we have
\begin{equation}\label{eq15}
B_p^{\mathbb{D}^*}(z) = |\log(|z|^2) |^p\beta_p^{\mathbb{D}^*}(z),
 \;\;\; \text{for} \;\; z\in \mathbb{D}^*,
\end{equation}
where
\begin{equation}\label{eq16}
\begin{split}
 \beta_p^{\mathbb{D}^*}(z) = \sum_{l=1}^{\infty}
 (c_l^{(p)})^2 |z|^{2l} , \;\;\;
 \text{with} \;\;\;
 c_l^{(p)} = \bigg( \frac{l^{p-1}}{2\pi(p-2)!} \bigg)^{1/2}.
\end{split}
\end{equation}
From (\ref{eq15}) and (\ref{eq16}), we know
$\displaystyle\lim_{z\rightarrow 0} B_p^{\D^*}(z) =0$ for any $p$.
Thus, $B_p^{\D^*}$ has no strictly positive lower bound.
 This is the main difficulty to establish Theorem 1.2.

By \cite[Proposition 3.3]{AMM21}, for any
$m\in \mathbb{N}, 0<b <1$ and $0<\gamma <\frac{1}{2}$, there exists
 $\varepsilon = \varepsilon(b,\gamma) >0$ such that
\begin{equation}\label{eq12}
\bigg \| B_p^{\mathbb{D}^*} (z)
- \frac{p-1}{2\pi} \bigg\|_{ C^m( \{be^{- p^{\gamma}}
\leq |z| <1\}, \omega_{\mathbb{D}^*} ) }
= \mathcal{O}(e^{-\varepsilon p^{1-2\gamma}}) \;\; \text{as} \;
p\rightarrow +\infty.
\end{equation}

After reducing to some $V_j$ and identifying the geometric data on
$\mathbb{D}^*_{4r}$ and $\Sigma$ via (\ref{eq8}),
by (\ref{eq1}) and (\ref{eq7}),
we have
\begin{equation}\label{eq11}
\begin{split}
& \frac{ J^*_{p,(2)} \omega_{FS,p}(z) } { p \;
\omega_{\mathbb{D}^*} (z)  }  - \frac{1}{2\pi}
  =  \frac{ i \partial\bar{\partial}
   \log (B^{\mathbb{D}^*}_p(z) )}{ 2\pi p \;
   \omega_{\mathbb{D}^*}(z) }  + \mathcal{O}(p^{-\infty}) \\
& = \frac{1}{2\pi p} |z|^2 \big( \log(|z|^2) \big)^2
\frac{ B_p^{\mathbb{D}^*}(z) \cdot
\frac{\partial^2}{\partial z \partial\bar{z}}
B^{\mathbb{D}^*}_p(z) - \frac{\partial}{\partial z}
B_p^{\mathbb{D}^*}(z) \cdot  \frac{\partial}{\partial \bar{z}}
B^{\mathbb{D}^*}_p(z) }{(B_p^{\mathbb{D}^*}(z))^2}
  + \mathcal{O}(p^{-\infty}) .
\end{split}
\end{equation}
As $|z|\log|z|^2 \frac{\partial}{\partial z}$ is an orthonormal frame of $(T^{(1,0)} \D^*, \omega_{\D^*})$,
from  (\ref{eq12}) with $m=1,2$, we have
\begin{equation}\label{eq003}
\begin{split}
&\sup_{  be^{- p^{\gamma}} \leq |z| <1 }  |z|^2
\big( \log(|z|^2) \big)^2 \bigg|\frac{\partial^2}
{\partial z \partial\bar{z}} B^{\mathbb{D}^*}_p(z)\bigg|
= \mathcal{O}(e^{-\varepsilon p^{1-2\gamma}}) \;\; \text{as} \;
 p\rightarrow +\infty. \\
&\sup_{ be^{- p^{\gamma}} \leq |z| <1 } |z|
   \big| \log(|z|^2) \big| \bigg|\frac{\partial}{\partial z}
   B_p^{\mathbb{D}^*}(z)\bigg|
= \mathcal{O}(e^{-\varepsilon p^{1-2\gamma}}) \;\; \text{as} \;
p\rightarrow +\infty.
\end{split}
\end{equation}
From  (\ref{eq12}) with $m=0$, (\ref{eq11}) and (\ref{eq003}),  we get
\begin{equation}\label{eq13}
\sup_{  be^{- p^{\gamma}} \leq |z| <1 }  \Bigg | \frac{ J^*_{p,(2)}
\omega_{FS,p}(z) } { p\; \omega_{\mathbb{D}^*} (z)  }
 - \frac{1}{2\pi} \Bigg | = \mathcal{O}(p^{-\infty}).
\end{equation}
Thus, in order to prove Theorem \ref{thm_MainThm} it suffices to show that
\begin{equation}\label{eq14}
\sup_{ 0 < |z| < be^{- p^{\gamma}} }
\Bigg | \frac{ J^*_{p,(2)} \omega_{FS,p}(z) }
{ p\; \omega_{\mathbb{D}^*} (z)  }   \Bigg | = \mathcal{O}(p^3 ).
\end{equation}
From (\ref{eq1}), (\ref{eq7}), (\ref{eq001}) and (\ref{eq15}),
we get
\begin{equation}\label{eq19}
\frac{ J^*_{p,(2)} \omega_{FS,p} (z) }{ p\;
\omega_{\mathbb{D}^*}(z) }
= |z|^2 \big( \log |z|^2 \big)^2 \cdot
 \frac{ \frac{\partial^2}{\partial z\partial\bar{z}}
 \beta_p^{\mathbb{D}^*}(z)  \cdot \beta_p^{\mathbb{D}^*} (z)
- \frac{\partial}{\partial z} \beta_p^{\mathbb{D}^*}(z)
\frac{\partial}{\partial \bar{z}} \beta_p^{\mathbb{D}^*}(z) }
{2\pi p\; \big(\beta_p^{\mathbb{D}^*}(z) \big)^2 }
+ \mathcal{O}(p^{-\infty}) .
\end{equation}

A key idea of \cite{AMM22} is that we use the precise formula
(\ref{eq16}), and only the lower degree (i.e., $l=\mathcal{O}(p)$) part
of (\ref{eq16}) plays the essential contribution in (\ref{eq19}).
We define $\delta_p$ as follows,
\begin{equation}\label{eq22}
\delta_p  = \bigg\lfloor  \frac{(p-2) }{ 2 |\log r|} \bigg \rfloor ,
\end{equation}
and $\lfloor x\rfloor$ is the integral part of $x\in \mathbb{R}$. Set
\begin{equation}\label{eq21}
\begin{split}
 & I_{1,p} : =  \sum_{l=1}^{\delta_p} \sum_{m=1}^{\delta_p} l(l-m)
  ( c_l^{(p)} )^2 (c_m^{(p)} )^2 |z|^{2l+2m},  \\
 & I_{2,p} : =  \sum_{l=1}^{\delta_p} \sum_{m=\delta_p +1}^{\infty}  l(l-m) ( c_l^{(p)} )^2 (c_m^{(p)} )^2 |z|^{2l+2m},  \\
 & I_{3,p} : =  \sum_{l=\delta_p +1}^{\infty} \sum_{m=1}^{\delta_p}  l(l-m) ( c_l^{(p)} )^2 (c_m^{(p)} )^2 |z|^{2l+2m},  \\
 & I_{4,p} : =  \sum_{l=\delta_p+1}^{\infty}
 \sum_{m=\delta_p+1}^{\infty}
 l(l-m) ( c_l^{(p)} )^2 (c_m^{(p)} )^2 |z|^{2l+2m} .
\end{split}
\end{equation}
From (\ref{eq16}) and (\ref{eq21}), we have
\begin{equation}\label{eq20}
\begin{split}
 & |z|^2 \bigg( \frac{\partial^2}{\partial z\partial\bar{z}}
 \beta_p^{\mathbb{D}^*}(z)  \cdot \beta_p^{\mathbb{D}^*} (z)
 - \frac{\partial}{\partial z} \beta_p^{\mathbb{D}^*}(z)
 \frac{\partial}{\partial \bar{z}} \beta_p^{\mathbb{D}^*}(z)
  \bigg) \\
 & = |z|^2 \sum_{l,m=1}^{\infty} l^2 (c_l^{(p)})^2
 (c_m^{(p)})^2 |z|^{2l+2m-2}
   - |z|^2 \sum_{l,m=1}^{\infty} lm( c_l^{(p)} )^2
   (c_m^{(p)} )^2 |z|^{2l+2m-2}  \\
 & = \sum_{l,m=1}^{\infty} l(l-m) ( c_l^{(p)} )^2
 (c_m^{(p)} )^2 |z|^{2l+2m} \\
 & = I_{1,p} + I_{2,p} + I_{3,p} + I_{4,p}.
\end{split}
\end{equation}
Thus, by (\ref{eq19}) and (\ref{eq20}) we get
\begin{equation}\label{eq24}
\frac{ J^*_{p,(2)} \omega_{FS,p} (z) }{ p \;
 \omega_{\mathbb{D}^*}(z) }
=  \big( \log |z|^2 \big)^2 \cdot  \frac{ I_{1,p}
+ I_{2,p} + I_{3,p} + I_{4,p} } { 2\pi p\;
\big(\beta_p^{\mathbb{D}^*}(z) \big)^2} +\mathcal{O}(p^{-\infty}).
\end{equation}
We will show the following result.

\begin{lem}\label{lem}
  (1) For $j=2,3,4$, we have
\begin{equation}\label{eq26}
\sup_{  0 < |z| < be^{- p^{\gamma}} }
\Bigg | \big( \log |z|^2 \big)^2 \cdot
\frac{I_{j,p} }{ \big( \beta_p^{\mathbb{D}^*}(z) \big)^2 }  \Bigg |
 = \mathcal{O}(p^{-\infty} ).
\end{equation}

(2) We have
\begin{subequations}
\begin{equation}\label{eq27}
\sup_{ 0 < |z| < 2 e^{-p} }  \Bigg | \big( \log |z|^2 \big)^2
\cdot \frac{I_{1,p} }{ \big( \beta_p^{\mathbb{D}^*}(z) \big)^2 }
 \Bigg | = \mathcal{O}(p^{-\infty} ),
\end{equation}
\begin{equation}\label{eq28}
\sup_{ e^{-p} < |z| < be^{- p^{\gamma}} }
\Bigg | \big( \log |z|^2 \big)^2 \cdot \frac{I_{1,p} }{ \big(
  \beta_p^{\mathbb{D}^*}(z) \big)^2 }  \Bigg | = \mathcal{O}(p^4 ).
\end{equation}
\end{subequations}
\end{lem}

\vskip 10pt
\noindent
{\bf How to get Theorem \ref{thm_MainThm} for Lemma \ref{lem}.}
Since
$\{z\in \mathbb{C}: 0<|z| < be^{-p^{\gamma}} \}
= \{z \in \mathbb{C}: 0<|z| < 2e^{-p}\} \cup \{z\in \mathbb{C}:
e^{-p} < |z| < be^{-p^{\gamma}} \}$, from (\ref{eq27})
and (\ref{eq28}), we get
\begin{equation}\label{eq006}
\sup_{ 0 < |z| < be^{- p^{\gamma}} }
\Bigg | \big( \log |z|^2 \big)^2 \cdot \frac{I_{1,p}}{
  \big( \beta_p^{\mathbb{D}^*}(z) \big)^2 }  \Bigg |
   = \mathcal{O}(p^4 ).
\end{equation}
From (\ref{eq24}), (\ref{eq26}) and (\ref{eq006}),
we get (\ref{eq14}).
The proof of Theorem \ref{thm_MainThm} is completed.

\subsection{Proof of Lemma \ref{lem}.}\label{s2.2}
\vskip 10pt
\noindent
From (\ref{eq21}), we have
\begin{equation}\label{eq29}
\begin{split}
 & | I_{1,p} | \leq \sum_{m=2}^{\delta_p} \big| 1-m \big|
 ( c_1^{(p)} )^2 (c_m^{(p)} )^2 |z|^{2+2m}
  + \sum_{l=2}^{\delta_p} \sum_{m=2}^{\delta_p }
  l\big|l-m \big| ( c_l^{(p)} )^2 (c_m^{(p)} )^2 |z|^{2l+2m}.
\end{split}
\end{equation}

We first estimate the first term on the right hand side of
(\ref{eq29}).
By (\ref{eq16}), we have
\begin{equation}\label{eq30}
  (c^{(p)}_1)^2 |z|^2 < \beta_p^{\mathbb{D}^*}(z),
\end{equation}
and
\begin{equation}\label{eq31}
  (c^{(p)}_m)^2 = \frac{ m^{p-1} }{ 2\pi(p-2)! }
  = \bigg( \frac{m}{m-1}\bigg)^{p-1} \cdot
  \frac{ (m-1)^{p-1} }{ 2\pi(p-2)! }
  = \bigg( \frac{m}{m-1}\bigg)^{p-1} (c^{(p)}_{m-1})^2 .
\end{equation}
From (\ref{eq16}) and (\ref{eq31}),
as $\frac{l}{l-1} \leq 2$ for $l\geq 2$, we have
\begin{equation}\label{eq0038}
\begin{split}
\sum_{l=2}^{\delta_p}   \big(c_l^{(p)} \big)^2  |z|^{2(l-1)}
& = \sum_{l=2}^{\delta_p}  \bigg( \frac{l}{l-1}\bigg)^{p-1}
\big(c_{l-1}^{(p)} \big)^2 |z|^{2(l-1)} \\
& \leq  2^{p-1}   \sum_{l=2}^{\delta_p}
 \big(c_{l-1}^{(p)} \big)^2 |z|^{2(l-1)}
\leq 2^{p-1} \beta_p^{\mathbb{D}^*} (z).
\end{split}
\end{equation}
From (\ref{eq22}), (\ref{eq30}) and (\ref{eq0038}), we have
\begin{equation}\label{eq38}
\begin{split}
\sum_{m=2}^{\delta_p} \big|1-m \big| ( c_1^{(p)} )^2
(c_m^{(p)} )^2 |z|^{2+2m}
&\leq (\delta_p-1) (c_1^{(p)})^2 |z|^2 \sum_{m=2}^{\delta_p}
(c_m^{(p)})^2 |z|^{2(m-1)} \cdot |z|^2  \\
& \leq (\delta_p -1) 2^{p-1} |z|^2  \; \beta_p^{\mathbb{D}^*}(z)^2 \\
&\leq p 2^{p-1} |z|^2 \beta_p^{\mathbb{D}^*}(z)^2 .
\end{split}
\end{equation}

\noindent
{\bf Proof of (\ref{eq27})}:
For $0<|z|<2 e^{-p}$, we have
\begin{equation}\label{eq35}
2^{p-1} |z| < 2^p e^{-p} = \bigg( \frac{2}{e} \bigg)^p.
\end{equation}
Thus, by (\ref{eq38}) and (\ref{eq35}), for $|z| < 2 e^{-p}$,
we get
\begin{equation}\label{eq36}
\begin{split}
 \sum_{m=2}^{\delta_p} \big| 1-m \big| ( c_1^{(p)} )^2
 (c_m^{(p)} )^2 |z|^{2+2m}
 \leq p \bigg( \frac{2}{e} \bigg)^p \; |z| \;
 \beta_p^{\mathbb{D}^*}(z)^2.
\end{split}
\end{equation}
Next we estimate the second term on the right hand side of
 (\ref{eq29}):
\begin{equation}\label{eq37}
\begin{split}
& \sum_{l=2}^{\delta_p} \sum_{m=2}^{\delta_p }  l\big|l-m\big|
 ( c_l^{(p)} )^2 (c_m^{(p)} )^2 |z|^{2l+2m} \\
& \leq \delta_p \cdot \delta_p \sum_{l=2}^{\delta_p}
\sum_{m=2}^{\delta_p} (c_l^{(p)} )^2 (c_m^{(p)})^2 |z|^{2l+2m}
\;\;\;\;\;\; (\text{as}\; |l-m|<\delta_p)\\
& = \delta_p^2 \; |z|^4 \; \sum_{l=2}^{\delta_p} (c_l^{(p)} )^2
 |z|^{2(l-1)} \cdot \sum_{m=2}^{\delta_p} (c_m^{(p)})^2
  |z|^{2(m-1)} .
\end{split}
\end{equation}
By (\ref{eq22}), (\ref{eq0038}), (\ref{eq35}) and (\ref{eq37}),
for $|z|<2e^{-p}$, as in (\ref{eq36}), we get
\begin{equation}\label{eq39}
\begin{split}
\sum_{l=2}^{\delta_p} \sum_{m=2}^{\delta_p } l\big|l-m\big|
 ( c_l^{(p)} )^2 (c_m^{(p)} )^2 |z|^{2l+2m}
& \leq  \delta_p^2 \; |z|^4 \cdot 2^{p-1}  \beta_p^{\mathbb{D}^*}
(z) \cdot 2^{p-1}  \beta_p^{\mathbb{D}^*} (z)  \\
& \leq  p^2 \cdot \bigg( \frac{2}{e} \bigg)^{2p} \; |z|^2 \;
 \beta_p^{\mathbb{D}^*}(z)^2 .
\end{split}
\end{equation}
Combining (\ref{eq29}), (\ref{eq36}) and (\ref{eq39}), we obtain
\begin{equation}\label{eq41}
\begin{split}
\Bigg| \big( \log |z|^2 \big)^2 \cdot \frac{I_{1,p}}{
  \big( \beta_p^{\mathbb{D}^*}(z) \big)^2 } \Bigg|
    \leq  p \bigg(\frac{2}{e}\bigg)^p \; |z| \big(\log |z|^2 \big)^2
+ p^2 \;\bigg(\frac{2}{e}\bigg)^{2p} \; |z|^2
\big(\log |z|^2 \big)^2 .
\end{split}
\end{equation}
Since the functions $f(x) = x(\log x^2)^2$ and $g(x)
= x^2(\log x^2)^2$ are bounded on the interval
$]0, 2 e^{-p}[ \subset ]0,e^{-2}[$, by (\ref{eq41}),
we get that there exists $C_1 >0$ such that
\begin{equation}\label{eq42}
\begin{split}
\sup_{  0<|z|< 2 e^{-p} } \Bigg| \big( \log |z|^2 \big)^2 \cdot
 \frac{I_{1,p}}{ \big( \beta_p^{\mathbb{D}^*}(z) \big)^2 } \Bigg|
 \leq C_1 \Bigg[ p\bigg(\frac{2}{e} \bigg)^p
 + p^2 \bigg(\frac{2}{e}\bigg)^{2p}  \Bigg] .
\end{split}
\end{equation}
By (\ref{eq42}), we get (\ref{eq27}).

\vskip 10pt
\noindent
{\bf Proof of (\ref{eq28}).}
From (\ref{eq16}), (\ref{eq22}) and (\ref{eq21}), we have
\begin{equation}\label{eq43}
\begin{split}
 |I_{1,p}| & \leq \sum_{l=1}^{\delta_p} \sum_{m=1}^{\delta_p }
 l\big|l-m \big| ( c_l^{(p)} )^2 (c_m^{(p)} )^2 |z|^{2l+2m}  \\
& \leq  \delta_p \cdot \delta_p
\sum_{l=1}^{\delta_p}
\sum_{m=1}^{\delta_p }  ( c_l^{(p)} )^2 (c_m^{(p)} )^2 |z|^{2l+2m}\\
& \leq  p^2 \;\beta_p^{\mathbb{D}^*}(z)^2.
\end{split}
\end{equation}
Therefore,
\begin{equation}\label{eq44}
\begin{split}
 \Bigg | \big( \log |z|^2 \big)^2 \cdot \frac{I_{1,p}}{
   \big( \beta_p^{\mathbb{D}^*}(z) \big)^2 }  \Bigg |
 \leq p^2 \big(\log|z|^2 \big)^2.
\end{split}
\end{equation}
For $z \in \{u \in \mathbb{C} : e^{-p} <|u|< b e^{-p^{\gamma}} \}$,
 we have $e^{-2p} < |z|^2 <1$. Thus, $-2p < \log |z|^2 <0$.
 Therefore, $\big(\log |z|^2 \big)^2 < 4p^2$.
Hence, by (\ref{eq44}), we have
\begin{equation}\label{eq45}
\begin{split}
 \sup_{ e^{-p} < |z| < b e^{-p^{\gamma}} } \Bigg |
 \big( \log |z|^2 \big)^2 \cdot \frac{I_{1,p}}{
  \big( \beta_p^{\mathbb{D}^*}(z) \big)^2 }  \Bigg |
 \leq p^2 \cdot 4p^2 = 4 p^4.
\end{split}
\end{equation}
From (\ref{eq45}), we get (\ref{eq28}).

\vskip 10pt
\noindent
{\bf Proof of Lemma \ref{lem}(1) for $j=2$.}
From (\ref{eq21}), we have
\begin{equation}\label{eq46}
\begin{split}
| I_{2,p} | &\leq    \sum_{l=1}^{\delta_p}
\sum_{m=\delta_p +1}^{\infty} l\big|l-m\big| ( c_l^{(p)} )^2
(c_m^{(p)} )^2 |z|^{2l+2m}  \\
& < \sum_{l=1}^{\delta_p} l (c_l^{(p)})^2 |z|^{2l} \cdot
\sum_{m=\delta_p +1}^{\infty} m (c_m^{(p)} )^2 |z|^{2m} ,
\end{split}
\end{equation}
as $l|l-m| =l(m-l) < lm$.

From  \cite[(3.70)]{AMM22}, there exists $C>0$ such that for
any $s\in\mathbb{N}^*$, $p\geq 2$, we have
\begin{equation}\label{eq48}
|c_s^{(p)} | \leq C p^{\frac{1}{2}} \frac{1}{(2r)^s}
\frac{1}{ \big( \log |2r|^2 \big)^{p/2} }.
\end{equation}

Next, we need the following

\begin{lem}\label{lem2}(\cite[(3.60)]{AMM22}).
For any $\tau \in \mathbb{N}$ fixed, we have
\begin{equation}\label{eq52}
\frac{1}{ \big| \log(|2r|^2) \big|^{p/2} }
\bigg( \frac{|z|}{ 2r }\bigg)^{\delta_p -\tau +1}
\leq C p^{-\frac{1}{2}} \frac{1}{2^{\alpha' p}}
\beta_p^{\mathbb{D}^*}(z)^{\frac{1}{2} }
\end{equation}
for all $p\gg 1$ and $|z| \leq c' p^{-A'}$, where $c'$ and
$A'$ are constants defined as follows:
\begin{equation}\label{eq-lem1}
  A' =\frac{1}{2\alpha'}, \;\;\; \alpha'
  =\frac{1}{4|\log r|} \;\;\; \text{and} \;\;\;
  c'= r e^{1/2\alpha'} \big| \log(|2r|^2) \big|^{1/2\alpha'} .
\end{equation}
\end{lem}

\noindent
\begin{proof} For the sake of completeness, we include
   a proof here.
By (\ref{eq22}), for any $\tau \in \mathbb{N}$ fixed, we obtain
\begin{equation}\label{eq-lem2}
  \alpha' p \leq \delta_p - \tau  \;\;\; \text{for} \; p \gg 1.
\end{equation}
Thus, by (\ref{eq-lem1}) and (\ref{eq-lem2}) for
 $\tau \in \mathbb{N}$ fixed, we have
\begin{equation}\label{eq-lem3}
  \bigg( \frac{|z|}{2r} \bigg)^{2(\delta_p -\tau)/p}
  \frac{1} {\big| \log(|2r|^2)\big|}
  \leq \bigg( \frac{|z|}{2r} \bigg)^{2\alpha'}
  \frac{1} {\big| \log(|2r|^2)\big|}
  \leq 2^{-2\alpha'} \frac{e}{p}
\end{equation}
for $p\gg 1$, $|z| \leq c' p^{-A'}$.
Recall that the Stirling formula states
\begin{equation}\label{eq-lem4}
  \frac{p^p}{p!} = (2\pi p)^{-1/2} e^p
  \big( 1 + \mathcal{O}(p^{-1}) \big) \;\;\;
   \text{as} \;\; p\rightarrow + \infty.
\end{equation}
By (\ref{eq16}), (\ref{eq-lem3}) and (\ref{eq-lem4}),
for any $\tau \in \mathbb{N}$ fixed, we have
\begin{equation}
  \begin{split}
    \frac{1}{ \big| \log(|2r|^2) \big|^{p/2} }
    \bigg( \frac{|z|}{ 2r }\bigg)^{\delta_p -\tau +1} &= \frac{1}{2r} \Bigg(  \bigg( \frac{|z|}{2r}
\bigg)^{2(\delta_p -\tau)/p} \frac{1} {\big|
     \log(|2r|^2)\big|} \Bigg)^{p/2} \big( 2\pi (p-2)!
     \big)^{1/2} c_1^{(p)}|z|  \\
    &\leq C p^{-\frac{1}{2}} \frac{1}{2^{\alpha' p}}
    \beta_p^{\mathbb{D}^*}(z)^{\frac{1}{2} }
  \end{split}
\end{equation}
for $p\gg 1$, $|z| \leq c' p^{-A'}$.
\end{proof}

\vskip 10pt
From now on, we assume $|z|\leq c' p^{-A'}$, $p \gg 1$.
By (\ref{eq16}) and (\ref{eq22}) we get
\begin{equation}\label{eq55}
\begin{split}
\sum_{l=1}^{\delta_p}  l  (c_l^{(p)})^2  |z|^{2l}
\leq   \delta_p \sum_{l=1}^{\delta_p} (c_l^{(p)})^2 |z|^{2l}
\leq   \delta_p \beta_p^{\mathbb{D}^*}(z)
\leq   p \beta_p^{\mathbb{D}^*}(z) .
\end{split}
\end{equation}
On the other hand, by (\ref{eq48}), we have
\begin{equation}\label{eq56}
\begin{split}
\sum_{m=\delta_p+1}^{\infty}  m (c_m^{(p)})^2 |z|^{2m}
& \leq \sum_{m=\delta_p +1}^{\infty}   m C p \frac{1}{(2r)^{2m}}
\frac{1}{ \big| \log(|2r|^2) \big|^p } |z|^{2m} \\
& =  C p  \frac{1}{ \big| \log(|2r|^2) \big|^p }
\sum_{m=\delta_p +1}^{\infty}  m \bigg( \frac{|z|}{ 2r }\bigg)^{2m} .
\end{split}
\end{equation}
Since for $0 \leq \xi <1$,
\begin{equation*}
\begin{split}
\sum_{q=N+1}^{\infty} q \xi^{q-1}
= \Bigg( \sum_{q=N+1}^{\infty} \xi^q \Bigg)'
 = \bigg( \frac{ \xi^{N+1} }{ 1 -\xi } \bigg)'
 = \frac{ (N+1) \xi^N - N \xi^{N+1}  }{ (1-\xi)^2  }
 \leq \frac{ (N+1) \xi^N } { (1-\xi)^2 },
\end{split}
\end{equation*}
we have
\begin{equation}\label{eq57}
\begin{split}
\sum_{m=\delta_p +1}^{\infty} m \bigg( \frac{|z|}{2r}\bigg)^{2m}
  &=  \bigg( \frac{|z|}{2r} \bigg)^2  \sum_{m=\delta_p +1}^{\infty}
  m \bigg( \frac{|z|}{ 2r }\bigg)^{2(m-1)}  \\
& \leq  \bigg( \frac{|z|}{ 2r }\bigg)^2  (\delta_p + 1)
  \frac{  \Big(\frac{|z|}{2r} \Big)^{2\delta_p} } {
    \Big( 1 - \big( \frac{|z|}{2r} \big)^2 \Big)^2 } \\
 &\leq     \bigg( \frac{|z|}{ 2r }\bigg)^2  (\delta_p + 1)
  \cdot \bigg( \frac{4}{3} \bigg)^2 \bigg(\frac{|z|}{2r}
   \bigg)^{2\delta_p} .
\end{split}
\end{equation}
Therefore, by (\ref{eq22}), (\ref{eq52}) for $\tau=1$, (\ref{eq56}) and (\ref{eq57}), we obtain
\begin{equation}\label{eq58}
\begin{split}
\sum_{m=\delta_p+1}^{\infty}  m (c_m^{(p)})^2 |z|^{2m}
& \leq   \bigg( \frac{4}{3} \bigg)^2  \bigg( \frac{|z|}{ 2r }\bigg)^2
 (\delta_p + 1)\cdot C p \frac{1}{ \big| \log(|2r|^2) \big|^p }
    \bigg(\frac{|z|}{2r} \bigg)^{2\delta_p}  \\
& \leq  \bigg( \frac{4}{3} \bigg)^2   \bigg( \frac{|z|}{ 2r }\bigg)^2
 (p + 1)  C'  \frac{1}{2^{2\alpha' p}} \beta_p^{\mathbb{D}^*}(z).
\end{split}
\end{equation}
From (\ref{eq46}), (\ref{eq55}) and (\ref{eq58}), 
 we get
\begin{equation}\label{eq61}
\begin{split}
 \Bigg | \big( \log |z|^2 \big)^2 \cdot \frac{ I_{2,p} }{
   \big( \beta_p^{\mathbb{D}^*}(z) \big)^2 }  \Bigg |
& \leq
\Big( \log |z|^2 \Big)^2    \bigg( \frac{|z|}{ 2r }\bigg)^2  \;
\bigg( \frac{4}{3} \bigg)^2  C'    p (p + 1)
\frac{1}{2^{2\alpha' p}} .
\end{split}
\end{equation}
Since the function $f(x) = x^2(\log x^2)^2$ is bounded on the
interval $]0,1[$, from (\ref{eq61}), we get
\begin{equation}\label{eq62}
\begin{split}
 \sup_{  0<|z| \leq c' p^{-A'}  } \Bigg | \big( \log |z|^2 \big)^2
  \cdot \frac{I_{2,p}}{ \big( \beta_p^{\mathbb{D}^*}(z) \big)^2 }
    \Bigg |
 \leq C'' p^2 \frac{1}{2^{2\alpha' p}} .
\end{split}
\end{equation}
Since $c' p^{-A'} > b e^{-p^{\gamma}}$ for
$p\gg 1$, we get (\ref{eq26}) for $j=2$ from (\ref{eq62}).

\vskip 10pt
\noindent
{\bf Proof of Lemma \ref{lem}(1) for $j=3$.}
We assume $|z|\leq c' p^{-A'}$, $p \gg 1$ as in Lemma \ref{lem2}.
From (\ref{eq21}), we have
\begin{equation}\label{eq63}
\begin{split}
|I_{3,p}|   & \leq   \sum_{l=\delta_p +1}^{\infty}
\sum_{m=1}^{\delta_p} l\big| l-m \big| ( c_l^{(p)} )^2
(c_m^{(p)} )^2 |z|^{2l+2m}   \\
& \leq \sum_{l=\delta_p + 1}^{\infty} l^2  (c_l^{(p)} )^2 |z|^{2l}
 \cdot \sum_{m=1}^{\delta_p}  (c_m^{(p)})^2 |z|^{2m},
\end{split}
\end{equation}
as $l|l-m| =l(l-m) <l^2$.

From (\ref{eq48}), we have
\begin{equation}\label{eq65}
\begin{split}
\sum_{l=\delta_p+1}^{\infty}  l^2 (c_l^{(p)})^2 |z|^{2l}
& \leq \sum_{l=\delta_p +1}^{\infty}   l^2 \; C p
\frac{1}{(2r)^{2l}} \frac{1}{ \big| \log(|2r|^2) \big|^p } \;
 |z|^{2l } \\
& =   C p  \frac{1}{ \big| \log(|2r|^2) \big|^p }
\sum_{l=\delta_p +1}^{\infty}  l^2 \bigg(
  \frac{|z|}{ 2r }\bigg)^{2l}  .
\end{split}
\end{equation}
By a simple computation, we get for $0\leq \xi <1$,
\begin{equation*}
\begin{split}
\sum_{q=N+1}^{\infty} q(q+1) \xi^{q-1}
&= \Bigg( \sum_{q=N+1}^{\infty} \xi^{q+1} \Bigg)'' \\
&= \bigg( \frac{ \xi^{N+2} }{ 1 -\xi } \bigg)'' \\
&=  \frac{ (N+2)(N+1) \xi^N - 2N(N+2)\xi^{N+1}
+ N(N+1) \xi^{N+2} }{ (1-\xi)^3 }   \\
& \leq \frac{(N+2) (N+1) \xi^N } { (1-\xi)^3 }.
\end{split}
\end{equation*}
Thus,
\begin{equation}\label{eq66}
\begin{split}
\sum_{l=\delta_p +1}^{\infty}l^2 \bigg( \frac{|z|}{ 2r }\bigg)^{2l}
& \leq  \bigg( \frac{|z|}{2r} \bigg)^2 \sum_{l=\delta_p +1}^{\infty}
 l(l+1) \bigg( \frac{|z|}{ 2r }\bigg)^{2(l-1)}  \\
& \leq  \bigg( \frac{|z|}{ 2r }\bigg)^2  (\delta_p+2) (\delta_p + 1)
  \frac{  \Big(\frac{|z|}{2r} \Big)^{2\delta_p} } {\Big( 1
  - \big( \frac{|z|}{2r} \big)^2 \Big)^3 }  \\
& \leq \bigg( \frac{|z|}{ 2r }\bigg)^2  (\delta_p + 2)
(\delta_p + 1)    \cdot   \bigg( \frac{4}{3} \bigg)^3
\bigg( \frac{|z|}{ 2r }\bigg)^{2 \delta_p }  .
\end{split}
\end{equation}
Therefore, by (\ref{eq65}) and (\ref{eq66}), we get
\begin{equation}\label{eq67}
\begin{split}
\sum_{l=\delta_p+1}^{\infty}  l^2 (c_l^{(p)})^2 |z|^{2l}
\leq   C  \bigg( \frac{|z|}{ 2r }\bigg)^2  (\delta_p + 2)
(\delta_p + 1) \bigg( \frac{4}{3} \bigg)^3   \cdot
p  \frac{1}{ \big| \log(|2r|^2) \big|^p }
\bigg( \frac{|z|}{ 2r }\bigg)^{2 \delta_p } .
\end{split}
\end{equation}
By (\ref{eq22}), (\ref{eq52}) for $\tau=1$ and (\ref{eq67}), we get
\begin{equation}\label{eq68}
\begin{split}
\sum_{l=\delta_p+1}^{\infty}  l^2 (c_l^{(p)})^2 |z|^{2l}
& \leq  C  \bigg( \frac{|z|}{ 2r }\bigg)^2  (\delta_p + 2)
(\delta_p + 1) \bigg( \frac{4}{3} \bigg)^3   \cdot C^2
\frac{1}{2^{2\alpha' p}}  \beta_p^{\mathbb{D}^*}(z)  \\
& \leq  C'  \bigg( \frac{|z|}{ 2r }\bigg)^2  (p+2) (p + 1)
 \frac{1}{2^{2\alpha' p}} \beta_p^{\mathbb{D}^*}(z) .
\end{split}
\end{equation}
From   (\ref{eq16}),
 (\ref{eq63}) and (\ref{eq68}), we get
\begin{equation}\label{eq70}
\begin{split}
\Bigg | \big( \log |z|^2 \big)^2 \cdot \frac{ I_{3,p} }{
  \big( \beta_p^{\mathbb{D}^*}(z) \big)^2 }  \Bigg |
\leq
\Big( \log |z|^2 \Big)^2    \bigg( \frac{|z|}{ 2r }\bigg)^2
\;\bigg[    C' (p+2)(p+1)   \frac{1}{2^{ 2\alpha' p} }  \bigg] .
\end{split}
\end{equation}
Since the function $f(x) = x^2(\log x^2)^2$ is bounded on the
interval $]0,1[$, from (\ref{eq70}), we get
\begin{equation}\label{eq71}
\begin{split}
& \sup_{ 0<|z| \leq c' p^{-A'} } \Bigg | \big( \log |z|^2 \big)^2
\cdot \frac{I_{3,p}}{ \big( \beta_p^{\mathbb{D}^*}(z) \big)^2 }
 \Bigg |
 \leq C_2 \; p^2 \frac{1}{2^{2\alpha' p}} .
\end{split}
\end{equation}
Since $c' p^{-A'} > b e^{-p^{\gamma}}$ for
$p\gg 1$, we get (\ref{eq26}) for $j=3$ from (\ref{eq71}).

\vskip 10pt
\noindent
{\bf Proof of Lemma \ref{lem}(1) for $j=4$.}
We assume $|z|\leq c' p^{-A'}$, $p \gg 1$ as in Lemma \ref{lem2}.
By (\ref{eq21}), we have
\begin{equation}\label{eq72}
\begin{split}
  | I_{4,p} |
& \leq  \sum_{l=\delta_p+1}^{\infty}  \sum_{m=\delta_p+1}^{\infty}
 l \big|l-m\big| ( c_l^{(p)} )^2 (c_m^{(p)} )^2 |z|^{2l+2m}  \\
& \leq \sum_{l=\delta_p + 1}^{\infty} l^2  (c_l^{(p)} )^2 |z|^{2l}
\cdot \sum_{m=\delta_p + 1 }^{ \infty}  (c_m^{(p)})^2 |z|^{2m}  \\
& \quad + \sum_{l=\delta_p +1}^{\infty} l (c_l^{(p)} )^2 |z|^{2l}
 \cdot \sum_{m=\delta_p + 1}^{\infty} m (c_m^{(p)} )^2 |z|^{2m} .
\end{split}
\end{equation}
From (\ref{eq48}),
as in (\ref{eq58}), we get
\begin{equation}\label{eq51}
\begin{split}
\sum_{m=\delta_p+1}^{\infty} (c_m^{(p)})^2 |z|^{2m}
& \leq \sum_{m=\delta_p +1}^{\infty} C p \frac{1}{(2r)^{2m}}
\frac{1}{ \big| \log(|2r|^2) \big|^p } |z|^{2m} \\
&\leq Cp \frac{1}{ \big| \log(|2r|^2) \big|^p }
4\cdot  \bigg( \frac{|z|}{ 2r }\bigg)^{2\delta_p +2} .
\end{split}
\end{equation}
From (\ref{eq52}) for $\tau =1$ and (\ref{eq51}), we get
\begin{equation}\label{eq53}
\begin{split}
\sum_{m=\delta_p+1}^{\infty} (c_m^{(p)})^2 |z|^{2m}
& \leq 4C \frac{1}{2^{ 2\alpha' p} } \beta_p^{\mathbb{D}^*}(z)
\bigg( \frac{|z|}{ 2r }\bigg)^2  .
\end{split}
\end{equation}
From (\ref{eq58}), (\ref{eq68}), (\ref{eq72}) and (\ref{eq53}), we get
\begin{equation}\label{eq73}
\begin{split}
|I_{4,p}| \leq  & C' \bigg(\frac{|z|}{2r} \bigg)^2 (p+2)(p+1)
\frac{1}{2^{2\alpha' p }} \beta_p^{\mathbb{D}^*}(z) \cdot
4C \frac{1}{2^{2\alpha' p} } \beta_p^{\mathbb{D}^*}(z)
 \bigg(\frac{|z|}{2r} \bigg)^2   \\
& + \Bigg[ \bigg(\frac{4}{3}\bigg)^2 C' \bigg( \frac{|z|}{2r} \bigg)^2
 (p + 1) \frac{1}{2^{2\alpha' p}}
 \beta_p^{\mathbb{D}^*}(z) \Bigg]^2  \\
\leq & C'' p^2 \frac{1}{2^{4\alpha' p}}
\bigg(\frac{|z|}{2r} \bigg)^4 \beta_p^{\mathbb{D}^*}(z)^2 .
\end{split}
\end{equation}
Thus, by (\ref{eq73}) we obtain
\begin{equation}\label{eq74}
\begin{split}
\Bigg | \big( \log |z|^2 \big)^2 \cdot \frac{I_{4,p} }{
  \big( \beta_p^{\mathbb{D}^*}(z) \big)^2 }  \Bigg |
& \leq C'' p^2 \frac{1}{2^{4\alpha' p}}
\bigg(\frac{|z|}{2r} \bigg)^4 \Big( \log |z|^2 \Big)^2 .
\end{split}
\end{equation}
Since the function $f(x) = x^4 (\log x^2)^2$ is bounded on
the interval $]0,1[$, from (\ref{eq74}), we get
\begin{equation}\label{eq75}
\begin{split}
\sup_{  0<|z| \leq c' p^{-A'} } \Bigg | \big( \log |z|^2 \big)^2
\cdot \frac{I_{4,p} }{ \big( \beta_p^{\mathbb{D}^*}(z) \big)^2 }
 \Bigg |
\leq C''' p^2 \frac{1}{2^{4\alpha' p}} .
\end{split}
\end{equation}
Since $c' p^{-A'} > b e^{-p^{\gamma}}$
 for $p\gg 1$, we get (\ref{eq26}) for $j=4$ from (\ref{eq75}).

The proof of Lemma \ref{lem} is completed.

\subsection{Proof of Corollary \ref{cor}}\label{s2.3}
Let $\Gamma$ be a geometrically finite Fuchsian group of the
 first kind, without elliptic elements,
then as explained in \cite[p956]{AMM21},
$\Sigma:= \Gamma \backslash \mathbb{H}$ can be compactified
by finitely many points $D=\{a_1,\cdots, a_N\}$
into a compact Riemann surface $\overline{\Sigma}$ with genus $g$
such that the following equivalent conditions (i)-(iv)
 are fulfilled:

\vskip 5pt
(i) $\Sigma = \overline{\Sigma} \setminus D$ admits a complete
K\"ahler\textendash Einstein metric $\omega_{\Sigma}$ with
$Ric_{\omega_{\Sigma}} = - \omega_{\Sigma}$,

(ii) $2g-2+N >0$,

(iii) the universal cover of $\Sigma$ is the upper-half plane
$\mathbb{H}$,

(iv) $L=K_{\overline{\Sigma}} \otimes
\mathcal{O}_{\overline{\Sigma}}(D)$ is ample,
 where $K_{\overline{\Sigma}} = T^{*(1,0)} \overline{\Sigma}$
  is the canonical line bundle on $\overline{\Sigma}$.

\vskip 5pt
\noindent
From \cite[Lemma 6.2]{AMM21}, there exists a singular Hermitian
metric on $L$ such that $(\Sigma, \omega_{\Sigma})$ and the formal
 square root of $(L,h)$ satisfy the condition \ref{alpha}
  and \ref{beta}.

From \cite[the proof of Lemma 6.2]{AMM21}, $\omega_{\Sigma}$ is
the K\"ahler\textendash Einstein form of constant negative curvature $-4$,
 induced by the Poincar\'e form
 $\omega_{\mathbb{H}} = \frac{ idz \wedge d\bar{z}}{ 4|{\rm Im} z|^2}$
 on $\mathbb{H}$,
and every $a\in D$ has a coordinate neighborhood
$(\overline{U}_a, z)$ in $\overline{\Sigma}$ such that
in this coordinate $\omega_{\Sigma}$ is exactly given by
 $\omega_{\mathbb{D}^*}(z)$ on
  $U_a = \overline{U}_a \setminus \{a\}$. Moreover, the curvature
  of the line bundle $(L|_{\Sigma}, h)$ is given by
  $-2i\omega_{\Sigma}$ and thus the curvature of the
  (formal) square root of $(L,h)$ is $-i\omega_{\Sigma}$.
Therefore, from \cite[Corollary 2.4]{AMM21}, we get that for any
 $k,m\in \mathbb{N}$ and any compact set $K\subset \Sigma$,
\begin{equation}\label{eq76}
B_p(x) = \frac{1}{\pi} p - \frac{1}{2\pi} + \mathcal{O}(p^{-k}), \;\;\;\;
\text{in} \;\; C^m(K) \;\; \text{as} \;\;
p\rightarrow +\infty.
\end{equation}
Hence, from (\ref{eq6}), we get
\begin{equation}\label{eq77}
\sup_{  K }  \Bigg | \frac{ J^*_{p,(2)} \omega_{FS,p}(z) } {
  2p \; \omega_{\Sigma} (z)  }  - \frac{1}{2\pi} \Bigg |
   = \mathcal{O}(p^{-\infty}) .
\end{equation}
On the coordinate neighborhood $U_a$ of $a\in D$, we have
$\omega_{\Sigma} = \omega_{\mathbb{D}^*}$.
Combine Theorem \ref{thm_MainThm} and (\ref{eq77}), we get
\begin{equation}\label{eq79}
\sup_{ \Sigma }  \Bigg | \frac{ J^*_{p,(2)} \omega_{FS,p}(z) } {
  p \; \omega_{\Sigma} (z)  } \Bigg | = \mathcal{O}(p^3).
\end{equation}

As explained in \cite[(6.36)]{AMM21},
 the space $\mathcal{S}_{2p}^{\Gamma}$ of cusp forms (Spitzenformen)
  of weight $2p$ of $\Gamma$ endowed with the Petersson scalar
  product is isometric with the space of ${\bL}^2$ holomorphic
  sections of $L^p$ over $\Sigma$ (cf. (1.2)) with respect to
  the volume form $\omega_{\Sigma}$ and the metric $h^p$ on $L^p$.
The Bergman kernel $B_p$ of $H^0_{(2)}(\Sigma, L^p)$ can
be identified to the Bergman kernel $S_p^{\Gamma}$
of $\mathcal{S}_{2p}^{\Gamma}$ (cf.  \cite[(6.37)]{AMM21}).
Thus, by (\ref{eq6}) and (\ref{eq.berg}), we get
\begin{equation}\label{eq.berg1}
\omega_{\Sigma}^{Ber,p} = \frac{i}{2\pi} \partial \bar{\partial} \log (B_p(z)),
\qquad 
\frac{1}{p} J^*_{p,(2)} \omega_{FS,p}
= \frac{1}{\pi} \omega_{\Sigma} + \frac{1}{ p} \omega_{\Sigma}^{Ber,p}.
\end{equation}
From (\ref{eq79}) and (\ref{eq.berg1}), we get (\ref{eq10}). 
The proof of Corollary \ref{cor} is completed.

\providecommand{\href}[2]{#2}

\end{document}